\documentclass[11pt]{article}

\usepackage{amsmath,amsthm,amsfonts,bbold}

\usepackage{amssymb}

\usepackage{epsfig}

\usepackage{rotating}

\usepackage{subfigure}

\newcommand{\fx} {\mathbf{x}}

\newcommand{\fa} {\mathbf{\alpha}}

\begin{document}

\title{Stochastic global optimization as a filtering problem}
\author{Panos Stinis \\ 
Department of Mathematics \\
University of Minnesota \\
    Minneapolis, MN 55455} 
\date {}

\maketitle

\begin{abstract}
We present a reformulation of stochastic global optimization as a filtering problem. The motivation behind this reformulation comes from the fact that for many optimization problems we cannot evaluate exactly the objective function to be optimized. Similarly, we may not be able to evaluate exactly the functions involved in iterative optimization algorithms. For example, we may only have access to noisy measurements of the functions or statistical estimates provided through Monte Carlo sampling. This makes iterative optimization algorithms behave like stochastic maps. Naive global optimization amounts to evolving a collection of realizations of this stochastic map and picking the realization with the best properties. This motivates the use of filtering techniques to allow focusing on realizations that are more promising than others. In particular, we present a filtering reformulation of global optimization in terms of a special case of sequential importance sampling methods called particle filters. The increasing popularity of particle filters is based on the simplicity of their implementation and their flexibility. For parametric exponential density estimation problems, we utilize the flexibility of particle filters to construct a stochastic global optimization algorithm which converges to the optimal solution appreciably faster than naive global optimization. Several examples are provided to demonstrate the efficiency of the approach.     
\end{abstract}


\section*{Introduction}
Optimization problems are ubiquitous in science and engineering ranging from portfolio optimization to space craft trajectory design to DNA studies and computer science \cite{horst, nocedal, papadimitriou}. As a result of the vast number of optimization applications and their related intricacies, the construction of optimization algorithms continues to be a subject of intense research. The main problem in optimization is to locate, usually through an iterative procedure, the optimal values (minima or maxima) of an objective function which depends on a number of parameters. The related problem of estimating the values of the parameters for which the objective function attains its optimal value is of equal importance in applications.  

For many optimization problems the function to be optimized is a function of several variables and can have multiple local extrema. Usually, optimization algorithms are guaranteed to find a local extremum when started from an arbitrary initial condition. In order to reach a global extremum we should use an optimization algorithm with multiple starting values and choose the best solution. This is called naive global optimization. In naive global optimization each solution starting from a different initial condition evolves independently of the others. Intuitively, we expect that a global optimization algorithm should benefit from the interaction of the different solutions. There is a category of global (deterministic and stochastic) optimization algorithms (e.g. branch and bound, tabu methods, genetic algorithms, ant colony optimization, swarm optimization, simulated annealing, hill climbing) \cite{eiben,horst,zhigljavsky} which allow for the different solutions to interact with the purpose of allocating more resources in areas of the parameter space which seem more promising in an optimization sense. 

In the current work, we propose an algorithm that belongs in this category. For many optimization problems we cannot evaluate exactly the objective function to be optimized. Similarly, we may not be able to evaluate exactly the functions involved in iterative optimization algorithms. For example, we may only have access to noisy measurements of the functions or statistical estimates provided through Monte Carlo sampling. This makes iterative optimization algorithms behave like stochastic maps. The corresponding global optimization problem also becomes stochastic. The algorithm we present here is based on the reformulation of the stochastic global optimization problem as a filtering problem. In particular, we present a filtering reformulation of stochastic global optimization in terms of a special case of sequential importance sampling methods called particle filters \cite{gordon,liu}. In this setting, the desired interaction of the solutions starting from different initial conditions is effected through the filtering step (see Section \ref{reformulation0} for more details).

The increasing popularity of particle filters is based on the simplicity of their implementation and their flexibility. The generic particle filter reformulation will not necessarily lead to an algorithm that converges faster than the naive global optimization algorithm. However, we are able to exploit the flexibility of the particle filter reformulation to construct a modified particle filter algorithm that converges appreciably faster than the naive global optimization algorithm. We provide several examples of parametric exponential density estimation of varying difficulty to demonstrate the efficiency of the approach.

The paper is organized as follows. Section \ref{local_global} discusses local and stochastic global optimization problems. In Section \ref{reformulation0} we present the reformulation of stochastic optimization problem as a filtering problem and in particular a reformulation which uses particle filters. In Section \ref{application0} we apply the particle filter reformulation to the problem of estimating the parameters of an exponential density. Section \ref{numerical0} contains numerical results for several examples to illustrate the efficiency of the proposed approach. Finally, in Section \ref{discussion} we provide a discussion and some directions for future work.

\section{Stochastic local and global optimization}\label{local_global}
Assume that we are given a scalar objective function $H(x)$ depending on $d$ parameters $x_1, \ldots, x_d.$ The purpose of the optimization problem is to compute the optimal (maximal or minimal) value of $H(x),$ say $\text{op}(H(x)),$ as well as the optimal parameter vector $\hat{x}=(\hat{x}_1,\ldots, \hat{x}_d)$ for which $H(\hat{x})=\text{op}(H(x)).$  For the sake of clarity and without loss of generality, let us assume that we are interested in obtaining the minimum value of the objective function $H(x).$

A generic optimization algorithm attempts to locate the optimal parameter vector in the following way:

\vskip14pt
{\bf Generic optimization algorithm}
\begin{enumerate}
\item
Pick an initial approximation $x^{(0)}.$ 
\item
For $k \ge 0,$ compute $x^{(k+1)}=f(x^{(k)}),$ where $f(x)=(f_1(x),\ldots,f_d(x))$ is a $d$-dimensional vector-valued function. Different optimization algorithms use a different function $f(x).$
\item
Evaluate $H(x^{(k+1)}).$ If $H(x^{(k+1)})$ satisfies a stopping criterion or $k+1$ is equal to a maximum number of iterations $k_{max}$, stop. Else, set $ k = k+1$ and proceed to Step 2.
\end{enumerate}
 
For many optimization problems, the function $H(x)$ has multiple local minima. The generic algorithm described above is, usually, able to locate one of the local minima, unless we are very lucky or know a lot about the problem at hand to guide the choice of $x^{(0)}.$  An obvious global version of the algorithm given above is as follows:

\vskip14pt
{\bf Naive global optimization algorithm}
\begin{enumerate}
\item
Pick a collection of initial approximations $x^{(0)}_1,\ldots, x^{(0)}_N.$

\item
Run, say for $M$ iterations, the generic optimization algorithm for the different initial conditions $x^{(0)}_1,\ldots, x^{(0)}_N.$ 
\item
Choose $\tilde{x}=\arg \underset{i=1\ldots N} {\min} H(x^{(M)}_i).$
\end{enumerate}
As stated in the introduction, for many optimization problems we only have access to a random estimate of the value of the objective function $H$ and of the optimization algorithm function $f.$ We will denote those as $H_s$ and $f_s$ respectively. The uncertainty can come from noisy measurements. Also, it can be due to Monte Carlo sampling error for cases where $H$ and/or $f$ involve expectations with respect to a probability density. Consequently, the iterative optimization sequence becomes stochastic. To be more precise, the optimization algorithm is modified as follows: 

\vskip14pt
{\bf  Stochastic optimization algorithm}
\begin{enumerate}
\item
Pick an initial approximation $x^{(0)}.$ 
\item
For $k \ge 0,$ compute $x^{(k+1)}=f_s(x^{(k)}).$
\item
Evaluate $H_s(x^{(k+1)}).$ If the value $H_s(x^{(k+1)})$ satisfies a stopping criterion or $k+1$ is equal to a maximum number of iterations $k_{max}$, stop. Else, set $ k = k+1$ and proceed to Step 2.
\end{enumerate}
Of course, there is an obvious stochastic version of the naive global optimization algorithm:
\vskip14pt
{\bf Naive stochastic global optimization algorithm}
\begin{enumerate}
\item
Pick a collection of initial approximations $x^{(0)}_1,\ldots, x^{(0)}_N.$
\item
Run $N$ realizations of the {\it stochastic} optimization algorithm, one for each of the different initial conditions $x^{(0)}_1,\ldots, x^{(0)}_N,$ say for $M$ iterations.
\item
Choose $\tilde{x}=\arg \underset{i=1\ldots N} {\min} H_s(x^{(M)}_i).$
\end{enumerate}

\section{Filtering reformulation of stochastic optimization}\label{reformulation0}
In the naive stochastic global optimization algorithm each realization of the stochastic optimization algorithm evolves independently of the other realizations. As it happens in many optimization problems, most realizations of a stochastic optimization algorithm will not start from an initial condition near a local minimum and thus will not contribute much to the exploration of minima in the space of parameter vectors. We would like to have a way to allocate more realizations in areas of the parameter vector space which seem more promising in an optimization sense. {\it We propose to do that by treating the stochastic global optimization problem as a filtering problem}. 

The filtering problem consists of the problem of incorporating the information from noisy observations of a system (noisy measurements of some function of the system) to correct the evolution trajectory of a system. The evolved system is usually stochastic and, in the case of discrete time formulations, a stochastic map. In the stochastic optimization algorithm context, we can think of: i) the iterative optimization algorithm, i.e. $x^{(k+1)}=f_s(x^{(k)}),$ as the stochastic map and  ii) the evaluation of the objective function, i.e. $H_s(x^{(k+1)}),$ as the noisy observation of the system. This allows us to recast (reformulate) the stochastic optimization problem as a filtering problem.

To make the reformulation more transparent we begin by recalling the abstract formulation of the filtering problem. Let $k$ be a discrete time index. Consider the following Markovian model which is motivated by the stochastic optimization algorithm:

\begin{align}
x^{(k+1)} &=f_s(x^{(k)}), \label{sys1} \\
y^{(k+1)}   &=H_s(x^{(k+1)})+v^{(k+1)}, \label{sys2}
\end{align}
where $v^{(k+1)}$ is a random variable with known properties. We have introduced the random variable $v^{(k+1)}$ to facilitate the formulation of the particle filter (see Section \ref{reformulation1}). Under appropriate regularity assumptions \cite{delmoral}, we can associate with \eqref{sys1} a state evolution density  $h(x^{(k+1)} | x^{(k)})$ and with \eqref{sys2} an observation density $g(y^{(k+1)} | x^{(k+1)}),$ i.e.
\begin{align}
x^{(k+1)} & \sim h(x^{(k+1)} | x^{(k)}),  \label{sys_den1} \\
y^{(k+1)} & \sim g(y^{(k+1)} | x^{(k+1)}). \label{sys_den2}
\end{align}  
Before we proceed with the details of the filtering problem we need to address the issue of the values of the observation sequence $y^{(k)}.$ In filtering problems, the observations are provided by some external measurement of the system. In the context of the stochastic optimization algorithm, there is no external measurement of the system. However, because we are dealing with an optimization problem (in our case minimization), the observation sequence $y^{(k)}$ is by construction a {\it non-increasing} sequence (in $k$). If we are looking for a maximum, then the sequence $y^{(k)}$ will be a non-decreasing sequence.

The filtering problem for the Markovian model \eqref{sys1}-\eqref{sys2} is to compute, for any time $t,$ the posterior distribution $p(x^{(0:t)} | y^{(1:t)})$ and/or the marginal distribution $p(x^{(t)} | y^{(1:t)}),$ where $x^{(0:t)}=(x^{(0)},\ldots,x^{(t)})$ and $y^{(1:t)}=(y^{(1)},\ldots,y^{(t)}).$ With the help of Bayes' theorem the posterior distribution $p(x^{(0:t)} | y^{(1:t)})$  can be written as (see e.g. \cite{doucet})
$$p(x^{(0:t)} | y^{(1:t)})= \frac{p(y^{(1:t)} | x^{(0:t)}) p(x^{(0:t)})}{\int p(y^{(1:t)} | x^{(0:t)}) p(x^{(0:t)})dx^{(0:t)}}.$$
The posterior distribution can be computed recursively by the formula 
$$ p(x^{(0:t+1)} | y^{(1:t+1)})= p(x^{(0:t)} | y^{(1:t)}) \frac{p(y^{(t+1)}| x^{(t+1)}) p(x^{(t+1)}|x^{(t)})}{p(y^{(t+1)}|y^{(1:t)})}.$$ 
The marginal distribution $p(x^{(t)} | y^{(1:t)})$ satisfies the recursion:
\begin{align}
Prediction: &\qquad  p(x^{(t)} | y^{(1:t-1)})= \int p(x^{(t)}|x^{(t-1)}) p(x^{(t-1)}|y^{(1:t-1)}) dx^{(t-1)},     \label{filter1} \\
Update:   & \qquad p(x^{(t)} | y^{(1:t)})=\frac{p(y^{(t)}| x^{(t)}) p(x^{(t)}|y^{(1:t-1)})}{\int p(y^{(t)} | x^{(t)}) p(x^{(t)}|y^{(1:t-1)})dx^{(t)}}.  \label{filter2}
\end{align}
The recursive formulas for the posterior and the marginal distribution can rarely be computed analytically, since for practical applications they involve the evaluation of complex high-dimensional integrals. Thus, we need to compute approximations to these distributions which should converge in some limit to the exact distributions. Sequential importance sampling (SIS) methods \cite{liu} have emerged as a flexible and powerful framework to construct such approximations. We should note here that SIS methods constitute a general framework for importance sampling, which is not restricted to filtering problems (for a nice presentation of the SIS framework in its full generality see Ch. 3 in \cite{liu}).   

We are now in a position to present the filtering reformulation of the stochastic optimization problem. It consists of {\it solving the filtering problem for 
\begin{align*}
x^{(k+1)} &=f_s(x^{(k)}), \\
y^{(k+1)}   &=H_s(x^{(k+1)})+v^{(k+1)}, 
\end{align*}
where $y^{(k+1)}$ is a {\it non-increasing} sequence.} We see that the stochastic optimization problem corresponds to a special case of the filtering problem since it includes the additional constraint of the non-increasing sequence of observations.

\subsection{Particle filter reformulation of stochastic global optimization}\label{reformulation1}

The approximation we will use to solve the filtering problem described above is a special case of the SIS formalism known as a particle filter (or bootstrap filter) \cite{gordon}. Assume (recall Eqs. \eqref{sys_den1}-\eqref{sys_den2}) that we are given a model for the state evolution density  $h(x^{(k+1)} | x^{(k)})$ and the observation density $g(y^{(k+1)} | x^{(k+1)}).$ In particular, let $x^{(k+1)}  \sim h(x^{(k+1)} | x^{(k)})$ and 
$y^{(k+1)}  \sim g(y^{(k+1)} | x^{(k+1)}).$ Suppose that at time $t$ we have a collection of $M$ random samples (particles) $(x^{(t)}_1,\ldots, x^{(t)}_M)$ which follow approximately the current marginal distribution $p(x^{(t)} | y^{(1:t)}).$ The paper \cite{gordon} suggests the following updating procedure after $y^{(t+1)}$ is observed:
\vskip14pt
{\bf Particle filter algorithm}
\begin{enumerate}
\item
Draw $x^{(t+1)}_{*j}$ from the state density $h(x^{(t+1)} | x^{(t)}_j),$ for $j=1,\ldots,M.$
\item
Assign to each draw the weight (also known as likelihood) $w^{*}_j= g(y^{(t+1)} | x^{(t+1)}_{*j}),$ for $j=1,\ldots,M.$
\item
Compute the normalized weights $w_j=\frac{w^{*}_j}{\sum_{l=1}^M w^{*}_l},$ for $j=1,\ldots,M.$
\item
Resample from $(x^{(t+1)}_{*1},\ldots, x^{(t+1)}_{*M})$ with probability proportional to $w_j, \; j=1,\ldots,M$ to produce new samples $(x^{(t+1)}_{1},\ldots, x^{(t+1)}_{M})$ at time $t+1.$
\item
Set $t=t+1$ and proceed to Step 1. 
\end{enumerate}
The particle filter approximates the marginal distribution $p(x^{(t)} | y^{(1:t)})$ by the distribution $\tilde{p}(x^{(t)} | y^{(1:t)})=\frac{1}{M}\sum_{i=1}^M \delta_{x^{(t)}_i}.$ It can be shown \cite{crisan}, under appropriate regularity assumptions on the state evolution density  $h(x^{(k+1)} | x^{(k)})$ and the observation density $g(y^{(k+1)} | x^{(k+1)}),$  that $\underset{M \rightarrow \infty}\lim \tilde{p}(x^{(t)} | y^{(1:t)})=p(x^{(t)} | y^{(1:t)})$ almost surely.

The particle filter construction appears as a natural candidate for the reformulation of the {\it global} stochastic optimization problem as a filtering problem. If one thinks of the particles in the particle filter construction as different initial conditions for the stochastic global optimization algorithm, then the particle filter becomes the filter analog of stochastic global optimization. Following the notation of the particle filter algorithm, we can specify the value of the observation at step $k$ as $y^{(k)}=\underset{j=1,\ldots,M} \min H_s(x^{(k)}_{*j}).$ With this choice for $y^{(k)},$ the particle filter formulation of the stochastic global optimization problem becomes

\vskip14pt
{\bf Particle filter algorithm for stochastic global optimization}
\begin{enumerate}
\item
Draw samples $x^{(0)}_{1},\ldots,x^{(0)}_{M}$ from an initial density $\mu_0(x).$ Set $t=0.$
\item
Draw samples $x^{(t+1)}_{*j}$ from the state density $h(x^{(t+1)}_j | x^{(t)}_j),$ for $j=1,\ldots,M.$
\item 
Compute $H_s(x^{(t+1)}_{*j})$ for $j=1,\ldots,M.$
\item 
Set $y^{(t+1)}=\underset{j=1,\ldots,M} \min H_s(x^{(t+1)}_{*j}).$
\item
Compute the weights $w^{*}_j(x^{(t+1)}_{*j})= g(y^{(t+1)} | x^{(t+1)}_{*j}),$ for $j=1,\ldots,M.$
\item
Compute the normalized weights $w_j(x^{(t+1)}_{*j})=\frac{w^{*}_j}{\sum_{l=1}^M w^{*}_l},$ for $j=1,\ldots,M.$
\item
Choose the estimate $\tilde{x}_{t+1}=\arg \underset{j=1\ldots M} {\max} w_j(x^{(t+1)}_{*j}).$ Equivalently, we can choose $\tilde{x}_{t+1}=\arg \underset{j=1,\ldots,M} \min H_s(x^{(t+1)}_{*j}).$
\item
If $t+1$ is equal to a maximum allowed number of iterations $t_{max}$ or $y^{(t+1)}$ satisfies a stopping criterion, terminate the algorithm. Else, proceed to next step.
\item
Resample from $(x^{(t+1)}_{*1},\ldots, x^{(t+1)}_{*M})$ with probability proportional to $w_j, \; j=1,\ldots,M$ to produce new samples $(x^{(t+1)}_{1},\ldots, x^{(t+1)}_{M})$ at time $t+1.$
\item
Set $t=t+1$ and proceed to Step 2. 
\end{enumerate}
The algorithm above computes a sequence of estimates $\tilde{x}_{t}, \; t=1,\ldots$ for the solution of the stochastic global optimization problem. Note that in Step 7, where we choose the estimate for the solution of the optimization problem, we provide two equivalent ways of determining the estimate. It is the same whether one chooses the particle with the smallest value of the objective function or the largest weight. This duality between minimizing the objective function and maximizing the weight (likelihood) of the observation of a filtering problem is not specific to this problem. The equivalence between least-squares problems and maximum likelihood estimation is well known in optimization literature (see e.g. \cite{nocedal}). Also, one can associate a filtering problem to a deterministic (Maslov) optimization process defined on a performance space  (see the appendix by Del Moral in \cite{maslov}). The duality between Maslov processes and filtering problems relies on idempotent analysis and, in particular the {\it Log-Exp} transform \cite{maslov}. The {\it Log-Exp} transform correspondence is present also in our particle filter reformulation of global optimization if we assume that the random variable $v^{(k+1)}$ in the observation process is Gaussianly distributed (see Section \ref{application2}).

The particle filter reformulation of global optimization does not suffer from any convergence drawbacks compared to the naive global optimization. To see this, we can assume for a moment that there is no observation step. Then, the particle filter algorithm reduces to naive global optimization. Naive global optimization will converge to a global minimum given an adequate number of particles. What the particle filter adds to the naive optimization is to pick, after every iteration, the solution that seems more promising in an optimization sense.   Since we have chosen the value of the observation to be equal to the {\it minimum} of the values of the objective function over all particles, the sequence of values picked for the objective function is guaranteed to be non-increasing. At worst, the value of the objective function will not change between iterations if the particles attempt moves in the parameter space that increase the value of the objective function. This can happen, especially for high-dimensional problems, where the local minima can be located in narrow valleys of the objective function \cite{nocedal}. The way we have defined the value of the observation ensures convergence to at least a local minimum. If, in addition, we have enough particles to explore in detail the parameter space, the particle filter will converge to the global minimum. However, the particle filter algorithm will not necessarily converge faster than the naive optimization algorithm. Yet, the reason we are presenting the particle filter reformulation of global optimization is because of its inherent flexibility. This means that starting from the generic form of the particle filter prescribed above, we can derive modified particle filters that exploit the specifics of the optimization algorithm used in the prediction step (e.g. Robbins-Monro, Kiefer-Wolfowitz, Levenberg-Marquardt algorithm etc.). This flexibility is exploited in Section \ref{application2.1} to construct a modified particle filter algorithm with faster convergence than the related naive optimization algorithm.

\section{Application to parametric exponential density estimation}\label{application0}
We apply the particle filter reformulation of stochastic global optimization to the problem of parametric density estimation for exponential densities. Exponential densities are, partly due to their nice mathematical features, widely used in the modeling of densities of systems of interacting variables in different contexts, ranging from Hamiltonian systems to image processing and bioinformatics (see \cite{jordan1} and references therein). As a result, there is an increased interest in algorithms for estimating and manipulating such densities numerically. In \cite{S05} we presented an algorithm  based on maximum likelihood for the estimation and renormalization (marginalization) of exponential densities. In the next section we present the main ideas of the work in \cite{S05}, in particular, how maximum likelihood estimation for exponential densities leads to an optimization problem.

\subsection{Parametric exponential density estimation as an optimization problem}\label{application1}
We begin our presentation with a few facts about families of exponential densities (for more details see e.g. \cite{brown,jordan1}). Let $\mathbf{x}=(x_1,\ldots,x_n)$ be an $n$-dimensional random vector taking values in $\mathbf{\Xi}^n \subseteq \mathbb{R}^n.$ The set $\mathbf{\Xi}^n=\Xi_1 \times \Xi_2 \times \ldots \Xi_n,$ where $x_1 \in \Xi_1,\ldots , x_n \in \Xi_n.$ Also, let $\psi_k(\mathbf{x}), k=1,\ldots,l$ be a collection of functions of $\mathbf{x}.$ The functions $\psi_k$ are known as potentials or sufficient statistics. Let $\psi=(\psi_1,\ldots,\psi_l)$ be the vector of potential functions. Associated with the vector $\psi$ is a vector $\mathbf{\alpha}=(\alpha_1,\ldots,\alpha_l)$ whose elements are called canonical or exponential parameters. The exponential family associated with $\mathbf{\psi}$ is the collection of density functions (parametrized by $\mathbf{\alpha}$) of the form
\begin{equation}
p(\mathbf{x},\mathbf{\alpha})=\frac{\exp(-\langle \mathbf{\alpha},\mathbf{\psi}(\mathbf{x})\rangle)}{Z(\mathbf{\alpha})},
\end{equation}
where $\langle \fa,\mathbf{\psi}(\fx) \rangle=\sum_{k=1}^l \alpha_k \psi_k(\fx)$  and $Z(\mathbf{\alpha})=\int_{\mathbf{\Xi}^n} \exp(-\langle \mathbf{\alpha},\mathbf{\psi}(\mathbf{x})\rangle)d\mathbf{x}.$ The exponential family is defined only for the set 
\begin{equation}
 \mathbf{A}= \{ \mathbf{\alpha} \in \mathbb{R}^l | Z(\mathbf{\alpha}) < \infty \}.
\end{equation}
Suppose that we are given a collection of $N$ independent samples of an $n$-dimensional random vector $\mathbf{x}$. In general, we do not know which density the samples are drawn from. There are many examples in practical applications where the random vector comes from a exponential density (see \cite{jordan1}). However, even if we do not know that the samples are drawn from an exponential density, we can try to fit an exponential density to the samples. This will become clearer when we formulate the moment-matching problem and exploit some properties of the exponential densities. The basic idea behind the algorithm we present here is to estimate the unknown parameter vector $\mathbf{\alpha}$ by maximizing the likelihood function of the samples. For a collection of $N$ independent samples of the random vector $\mathbf{x}$, the likelihood function $L$ is defined as (see e.g. \cite{lehman}) 
\begin{equation}
L=\prod_{j=1}^{N} p(\mathbf{x}_j,\alpha),
\end{equation}
where  $p(\mathbf{x}_j,\alpha)$ is the unknown exponential density whose parameters $\mathbf{\alpha}$ we wish to determine. We associate a potential function $\psi_k, k=1,\ldots,l$ with every parameter $\alpha_k.$ Maximization of $L$ with respect to the parameters $\alpha_k$ produces an estimate $\bar{\mathbf{\alpha}}$ for $\mathbf{\alpha}.$  Under suitable regularity conditions, the sequence of estimates $\bar{\mathbf{\alpha}}$ for increasing values of $N$ is asymptotically efficient and tends, with probability one, to a local maximum in parameter space. From now on we will use the notation $\mathbf{\alpha}$ instead of $\bar{\mathbf{\alpha}}$ to denote the maximum likelihood estimate of the parameters keeping in mind that this is only an estimate of the parameters. In addition, we will be working with the logarithm of the likelihood $\log L$, since it does not alter the position of the maximum and also leads to formulas that are more easily manipulated. Differentiation of $\log L$ with respect to the $\alpha_k$ and setting the derivative equal to zero results in 
\begin{equation}
\label{mom1}
E_{\mathbf{\alpha}}[\psi_k(\mathbf{x})]=\frac{1}{N} \sum_{j=1}^{N} \psi_k(\mathbf{x}_j), \quad k=1,\ldots,l
\end{equation}
where
\begin{equation}
E_{\mathbf{\alpha}}[\psi_k(\mathbf{x})]= \frac{\int_{\mathbf{\Xi}^n} \psi_k(\fx) \exp(-\langle \mathbf{\alpha},\mathbf{\psi}(\mathbf{x})\rangle) d\mathbf{x}}{\int_{\mathbf{\Xi}^n} \exp(-\langle \mathbf{\alpha},\mathbf{\psi}(\mathbf{x})\rangle) d\mathbf{x}} 
\end{equation}
is the expectation value of the function $\psi_k$ with respect to the density $p(\mathbf{x},\mathbf{\alpha}).$ The right side of (\ref{mom1}) is the average (moment) of the function $\psi_k$ as calculated from the given samples. The $l$ equations in (\ref{mom1}) define the {\it moment-matching} problem. What we want to do is to estimate the parameters $\mathbf{\alpha}$ so that the conditions in (\ref{mom1}) are satisfied. The question of whether such a problem has a solution, and if it does whether it is unique can be addressed through convex analysis (see \cite{S05} for more details). 

Now that we have defined the moment-matching problem we have to find a way to actually estimate the parameter vector $\fa.$ The equations (\ref{mom1}) contain, in general, nonlinear functions of the parameters. Moreover, except for very special cases, these nonlinear functions are unknown or very difficult to manipulate analytically. Thus, we have to tackle the problem of estimating the parameter vector numerically. We can define the $l$-dimensional vector  $\mathbf{f}(\fa)=(f_1(\fa),f_2(\fa),\ldots,f_l(\fa))$ as
\begin{equation}
\label{mom2}
f_k(\fa)=E_{\mathbf{\alpha}}[\psi_k(\mathbf{x})]-\frac{1}{N} \sum_{j=1}^{N} \psi_k(\mathbf{x}_j), \quad k=1,\ldots,l
\end{equation}
The moment-matching problem amounts to solving the system of (nonlinear) equations $f_k(\fa)=0, k=1,\ldots,l.$

\subsubsection{The Levenberg-Marquardt algorithm}\label{application1.1}
Two popular candidates to solve the system of nonlinear equations $f_k(\fa)=0, k=1,\ldots,l,$ are the method of steepest descent and Newton's method. However, both have their drawbacks. The method of steepest descent converges but can have very slow convergence, while Newton's method converges quadratically but it diverges if the initial guess of the solution is not good. We choose to solve the moment-matching problem as an optimization problem using the Levenberg-Marquardt (LM) algorithm (see e.g. \cite{nocedal}). This is a powerful iterative optimization algorithm that combines the advantages of the method of steepest descent and Newton's method. First, let us write the moment-matching problem as an optimization problem. Define the error function $\epsilon(\fa)$ as
\begin{equation}
\label{error}
\epsilon(\fa)=\frac{1}{2}\sum_{k=1}^l \epsilon_k^2=\frac{1}{2}\sum_{k=1}^l f_k^2(\fa),
\end{equation}
where $\epsilon_k=f_k(\alpha).$
The problem of minimizing $\epsilon(\fa)$ is equivalent to solving the system of equations $f_k(\fa)=0, k=1,\ldots,l$ i.e. the zeros of $\epsilon$ are solutions of the system $f_k(\fa)=0, k=1,\ldots,l$ and vice versa. The LM algorithm uses a positive parameter $\lambda$ to control convergence and the updates of the parameters at step $m+1$ are calculated through the formula
\begin{equation}
\label{lm1}
[J^T J+ \lambda  diag (J^T J) ](\alpha^{(m+1)}-{\alpha^{(m)}} )=-J^T  \mathbf{f}(\fa^{(m)}),
\end{equation}
where $J=\frac{\partial{f_i}}{\partial {\alpha_j}}|_{\fa=\fa^{(m)}}, i,j=1,\ldots,l$ is the Jacobian of $\mathbf{f}(\fa^{(m)})$ and $J^T$ its transpose. The matrix $diag (J^T J)$ is a diagonal matrix whose diagonal elements are the diagonal elements of $ (J^T J).$ In the literature, the name Levenberg-Marquardt is also used to denote the algorithm in (\ref{lm1}) with $diag (J^T J)$ replaced by the unit matrix $I.$ For the case where we use the unit matrix $I$ instead of $diag (J^T J),$ it is straightforward to see the connection with the methods of steepest descent and Newton's. For $\lambda=0$ the algorithm reduces to Newton's method, while for very large $\lambda$ we recover the steepest descent method. The modification (due to Marquardt) of using $diag (J^T J)$ becomes important in the case where $\lambda$ is large. In this case if we only used the unit matrix $I$ almost  all information coming from $ (J^T J)$ is lost. On the other hand, since $ (J^T J)$ provides information about the curvature of $\epsilon$, use of the matrix $ diag (J^T J)$ allows us to incorporate information about the curvature even in cases with large $\lambda.$ 

We have to prescribe a way of computing the Jacobian $J(\fa^{(m)})$. The element $J_{ij}$ of the  Jacobian is given by 
\begin{equation}
J_{ij}(\fa^{(m)})=- (E_{\fa^{(m)}}[\psi_i(\fx)\psi_j(\fx)]-E_{\fa^{(m)}}[\psi_i(\fx)]E_{\fa^{(m)}}[\psi_j(\fx)] )
\end{equation}
for $i,j=1,\ldots,l$ (note that the Jacobian is symmetric) So, all the quantities involved in equation (\ref{lm1}) can be expressed as expectation values with respect to the $m$-th step parameter estimate $\fa^{(m)}.$ More details about the implementation of the LM algorithm can be found in \cite{nocedal}.

If one can compute only Monte Carlo estimates of the necessary expectation values appearing in the LM algorithm, the question of convergence of the algorithm becomes more involved. Starting with the work of Robbins-Monro and Kiefer-Wolfowitz \cite{kushner}, there has been a vast literature on the subject of constructing convergent stochastic algorithms to obtain the zeros of a function for which we only have noisy observations. The optimization problem we are interested in fits in this category since we are looking for the zeros of the error function $\epsilon(\fa)$ for which we can only compute Monte Carlo estimates. It is not generally true that replacing the expectation values appearing in LM (or Newton's method) with Monte Carlo estimates will lead to a stochastic algorithm that converges to the true zero of $\epsilon(\fa)$ (see e.g. the discussion in \cite{gelman}). Convergence can be achieved if the number of Monte Carlo samples is increasing with iterations. In the limit of infinite number of samples, one will recover the deterministic LM algorithm. In particular, one can conceive of a scheme to increase the number of Monte Carlo samples whenever the error value becomes of the order of magnitude of the Monte Carlo error. An alternative scheme to increase the number of Monte Carlo samples was proposed in \cite{gelman}. On the other hand, if one can afford to have a large number of Monte Carlo samples for all iterations, the LM algorithm will reduce the error down to a value that is for all practical purposes zero (since the associated Monte Carlo error will be practically zero). We are not advocating here the use of the stochastic version of the LM algorithm for general problems. But, given the superiority of the deterministic LM algorithm for problems of moderate dimension like the ones we will be examining in Section \ref{numerical0} and the fact that we can afford a large number of Monte Carlo samples for {\it all} iterations, makes the stochastic LM algorithm a reasonable choice. We should also mention that what we refer to as a stochastic LM algorithm is different from the stochastic LM algorithm proposed by LeCun \cite{lecun-98b}.

\subsection{Reformulation of the density estimation optimization problem as a filtering problem}\label{application2}
For most cases, the expectation values involved in maximum likelihood estimation can only be computed through Monte Carlo sampling. This renders random both the update equation \eqref{lm1} for the coefficients and the equation \eqref{error} for the error. In addition, as mentioned above, maximum likelihood estimation is guaranteed to yield a local maximum in parameter space. Finding a global maximum likelihood estimate is equivalent to finding a global minimum, zero in this case, for the error function $\epsilon(\fa).$ This requires running the LM algorithm with different initial conditions and picking the best solution, i.e. the solution with the minimum value for the error function $\epsilon(\fa).$ The combination of randomness in the update and error equations and the need to find a global maximum of the likelihood (in the parameter space), means that the problem of global maximum likelihood estimation is equivalent to a stochastic global optimization problem. We will apply the filtering reformulation of the previous section in order to solve this stochastic global optimization problem. 

In order to reformulate the stochastic global optimization problem as a particle filter (see reformulation at the end of Section \ref{reformulation1}) we need to specify the state density $h(x^{(t+1)}_j | x^{(t)}_j),$ the objective function $H_s(x^{(t+1)}_{*j}),$ as well as the observation density $g(y^{(t+1)} | x^{(t+1)}_{*j}),$ where $j=1,\ldots,M.$  To conform  with the notation used above for exponential densities, this means that we need to specify $h(\fa^{(t+1)}_j | \fa^{(t)}_j),$ $H_s(\fa^{(t+1)}_{*j})$ and $g(y^{(t+1)} | \fa^{(t+1)}_{*j}).$

The state density  $h(\fa^{(t+1)}_j | \fa^{(t)}_j)$ is defined implicitly through the update equation \eqref{lm1}. Since the expectation values appearing in  \eqref{lm1} are estimated through Monte Carlo sampling,  Eq.\eqref{lm1} becomes random. In particular, the Jacobian $J_j(\fa_j^{(t)})$ for the $j$-th particle at iteration $t$ is replaced by its Monte Carlo estimate, say with $N_E$ samples, $\tilde{J}_j(\fa_j^{(t)}).$ Similarly, the moment matching (error) vector $\mathbf{f}(\fa_j^{(t)})$ is replaced by its Monte Carlo estimate $\tilde{\mathbf{f}}(\fa_j^{(t)}).$ For both Monte Carlo estimates the associated sampling error is $O(1/\sqrt{N_E})$ \cite{liu}.

For the objective function $H_s(\fa^{(t+1)}_{*j})$ we define $H_s(\fa^{(t+1)}_{*j})=(\frac{2*\tilde{\epsilon}(\fa_{*j}^{(t+1)})}{l})^{1/2},$ where $\tilde{\epsilon}(\fa_{*j}^{(t+1)})$ is the Monte Carlo estimate of the error ${\epsilon}(\fa_{*j}^{(t+1)})$ of the $j$-th particle. The quantity $(\frac{2*\tilde{\epsilon}(\fa_{*j}^{(t+1)})}{l})^{1/2}$ is the average error per moment for the $j$-th particle. The sampling error for $(\frac{2*\tilde{\epsilon}(\fa_{*j}^{(t+1)})}{l})^{1/2}$  is also $O(1/\sqrt{N_E}).$

We have to specify also the observation density $g(y^{(t+1)} | \fa^{(t+1)}_{*j}).$ Recall that we have defined the observation process $y^{(k+1)}$ as $y^{(k+1)}  =H_s(\fa^{(k+1)}_{*})+v^{(k+1)},$ where $\fa^{(k+1)}_{*}$ is the predicted state before the observation. We pick the random variable $v^{(k+1)}$ to be distributed as a Gaussian variable of mean zero and variance $O(1/N_E).$ This means that we can write the conditional observation density value for the $j$-th particle as $$g(y^{(t+1)} | \fa^{(t+1)}_{*j})=\frac{1}{\sqrt{2\pi \sigma^2}}\exp(-\frac{(y^{(t+1)} - H_s(\fa^{(t+1)}_{*j}))^2}{2*\sigma^2}),$$ where $\sigma^2=\frac{1}{N_E}.$ Also, as discussed at the end of Section \ref{reformulation1}, we pick the value of the observation $y^{(t+1)}=\underset{j=1,\ldots,M} \min H_s(\fa^{(t+1)}_{*j}).$ 

With the above choices, the particle filter algorithm becomes:

\vskip14pt
{\bf Particle filter algorithm for stochastic global optimization with the LM method}
\begin{enumerate}
\item
Draw samples (particles) $\fa^{(0)}_{1},\ldots,\fa^{(0)}_{M}$ from an initial density $\mu_0(x).$ Set $\lambda^{(0)}=\lambda_0$ and $t=0.$
\item
Use $[\tilde{J_j}^T \tilde{J_j}+ \lambda^{(t)}_j  diag (\tilde{J_j}^T \tilde{J_j}) ](\alpha_{*j}^{(t+1)}-{\alpha_{j}^{(t)}})= - \tilde{J_j}^T  \tilde{\mathbf{f}}(\fa_j^{(t)})$ to compute samples $\fa^{(t+1)}_{*j}$  for $j=1,\ldots,M.$  
\item 
Compute $H_s(\fa^{(t+1)}_{*j})=(\frac{2*\tilde{\epsilon}(\fa_{*j}^{(t+1)})}{l})^{1/2}$ for $j=1,\ldots,M.$ Determine $\lambda^{(t+1)}_j.$
\item 
Set $y^{(t+1)}=\underset{j=1,\ldots,M} \min H_s(\fa^{(t+1)}_{*j}).$
\item
Compute the weights $w^{*}_j(\fa^{(t+1)}_{*j})=\frac{1}{\sqrt{2\pi \sigma^2}}\exp(-\frac{(y^{(t+1)} - H_s(\fa^{(t+1)}_{*j}))^2}{2*\sigma^2}),$ where $\sigma^2=\frac{1}{N_E}$ for $j=1,\ldots,M.$
\item
Compute the normalized weights $w_j(\fa^{(t+1)}_{*j})=\frac{w^{*}_j}{\sum_{l=1}^M w^{*}_l},$ for $j=1,\ldots,M.$
\item
Choose the estimate $\tilde{x}_{t+1}=\arg \underset{j=1\ldots M} {\max} w_j(\fa^{(t+1)}_{*j}).$ Equivalently, we can choose $\tilde{x}_{t+1}=\arg \underset{j=1,\ldots,M} \min H_s(\fa^{(t+1)}_{*j}).$
\item
If $t+1$ is equal to a maximum allowed number of iterations $t_{max}$ or $y^{(t+1)}$ satisfies a stopping criterion, terminate the algorithm. Else, proceed to next step.
\item
Resample from $(\fa^{(t+1)}_{*1},\ldots, \fa^{(t+1)}_{*M})$ with probability proportional to $w_j, \; j=1,\ldots,M$ to produce new samples $(\fa^{(t+1)}_{1},\ldots, \fa^{(t+1)}_{M})$ at time $t+1.$
\item
Set $t=t+1$ and proceed to Step 2. 
\end{enumerate}

The algorithm described above allows one to allocate more particles in areas of the parameter space which seem promising in an optimization sense. A more careful inspection of the way the algorithm works reveals that there are two issues that can prevent the algorithm above to improve on the naive algorithm. 

First, the particle filter algorithm can converge fast to a local minimum that is not the best minimum that the optimization algorithm could have reached starting from a set of initial conditions. This phenomenon, called premature convergence, is well-known in the stochastic optimization literature, especially in the context of genetic algorithms \cite{eiben}. To be more precise, it is possible that early on in the iteration process one particle can dominate in the resampling step, because it appears to perform better initially. However, there may be other particles that would have performed better if we allowed them to evolve for more iterations. Such particles can vanish during the first few resampling steps because they have low weights. One way to avoid premature convergence in the context of the particle filter is to allow a few iterations without enforcing the filtering and resampling step so that the different particles can realize their potential to reach a minimum. After those iterations, the filtering and resampling steps can allocate more particles in the areas of the parameter space which are more promising.

Second, in the particle filter algorithm above, all the offspring particles are assigned the value of the LM parameter $\lambda$ of their parent. This means that all the offspring particles will evolve to the same new parameter vector. The only difference between the offspring particles will be with regards to their error value at the next step. This difference comes from the randomness in the Monte Carlo sampling. This may not be enough to give the particle filter an appreciable advantage as far as the convergence speed is concerned.     

In our numerical studies, which will be presented in Section \ref{numerical0}, we did not encounter the problem of premature convergence. However, as will be shown, the particle filter in its generic form did {\it not} achieve consistently a speedup of the convergence compared to the naive global optimization algorithm. Recall that the main advantage of the particle filter algorithm is the flexibility provided by the filtering and resampling steps. We exploit this flexibility to construct a modified particle filter algorithm which can increase appreciably the speed of convergence.  In some cases, the modified particle filter also achieved a smaller value of the error (compared to the naive algorithm) for a fixed maximum number of iterations.

\subsubsection{Modified particle filter algorithm}\label{application2.1} 
To motivate the modified particle filter algorithm we have to examine more carefully the way that the LM algorithm works. As mentioned in Section \ref{application1.1}, the LM algorithm is a hybrid of the steepest descent method and Newton's method. In particular, depending on the value of the parameter $\lambda,$ the LM algorithm can be brought closer to steepest descent or Newton's method. For large values of $\lambda$ it is close to the steepest descent method, while for $\lambda=0$ it reduces to Newton's method (assuming that the Jacobian is not singular). While the steepest descent method has guaranteed convergence the rate of convergence can be slow. On the other hand, Newton's method has a higher speed of convergence but whether it will converge or not is sensitive on the choice of initial conditions. 

The resampling step of the particle filter produces more copies (offspring) of a good, in an optimization sense, particle. The offspring particles inherit {\it all} the properties of the parent particle and, in particular, the same value for the parameter $\lambda,$ say $\lambda_P.$ However, nothing prevents us from assigning {\it different} values of $\lambda$ to the different offspring of the same particle. After all, $\lambda$ is an adaptive parameter of the algorithm which should be chosen so as to accelerate convergence. The resampling step of the particle filter allows us to {\it batch} the offspring of a particle according to their ancestry. Then, we can assign different values to the offspring within a batch. Let $(\fa^{(t+1)}_{1},\ldots, \fa^{(t+1)}_{B})$ be the batch of offspring of a good particle. We can assign to the $i$-th element of this batch value  $\lambda^{(t+1)}_{i},$ where $i=1,\ldots,B.$ A value $\lambda^{(t+1)}_{i} > \lambda_P$ means that the LM algorithm for the $i$-th offspring particle will behave closer to the steepest descent method than the parent particle. Similarly, a value $\lambda^{(t+1)}_{i} < \lambda_P$ means that the LM algorithm for the $i$-th offspring particle will behave closer to Newton's method than the parent particle. Since our goal is to accelerate convergence, i.e. bring the LM algorithm closer to Newton's method, we chose to assign the values of $\lambda$ to the offspring within a batch in the following manner: $\lambda^{(t+1)}_{i}= \lambda_P/\gamma^{i-1},$ with $\gamma > 1.$ This means that $\lambda^{(t+1)}_{1}=\lambda_P$ and $\lambda^{(t+1)}_{B}= \lambda_P/\gamma^{B-1}.$ For batches containing many offspring (large $B$), this procedure allows us to explore fast the neighborhood of a good particle. Note that such an exploration is impossible in the context of the naive global optimization algorithm without increasing tremendously the cost of the algorithm. On the other hand, the cost of this exploration is negligible in the particle filter framework. A good choice for $\gamma$ can be found based on knowledge of how the parameter $\lambda$ influences the behavior of the LM algorithm, in particular the determinant and condition number of the matrix $\tilde{J_j}^T \tilde{J_j}+ \lambda^{(t)}_j  diag (\tilde{J_j}^T \tilde{J_j}).$  The quality of the results appears to be pretty robust to the value of $\gamma.$ Also, it is encouraging that for all the numerical examples we could use the same value of $\gamma$ and obtain  equally good results.

\section{Numerical results}\label{numerical0}
We have applied the particle filter algorithm (both in its generic and modified forms) to four examples of varying difficulty. The first two examples involve the estimation of the parameters of a known exponential density. The last two involve the estimation of the  parameters of an exponential density so as to match certain moments of an unknown density. 

For all four examples the exponential density whose parameters are to be estimated is given by 
\begin{equation*}
p(\mathbf{x},\mathbf{\alpha})=\frac{\exp(-\langle \mathbf{\alpha},\mathbf{\psi}(\mathbf{x})\rangle)}{Z(\mathbf{\alpha})},
\end{equation*}
where $\fa=(\fa_1,\ldots,\fa_{24})$ and $\mathbf{x}=(x_1,\ldots,x_4).$ The potential function vector  $\mathbf{\psi}(\mathbf{x})$ is defined as 
\begin{equation}
\label{pots}
\begin{split}
\psi_{1-4}(\mathbf{x})=&  x_i, \; \text {for} \; i\,=1,\ldots,4 \\
\psi_{5-14}(\mathbf{x})=& x_i x_j, \; \text{for}\; i,j=1,\ldots,4 \; \text{and} \; j \geq i \\
\psi_{15-24}(\mathbf{x})=& x_i^2 x_j^2, \; \text{for}\; i,j=1,\ldots,4 \; \text{and} \; j \geq i 
\end{split}
\end{equation}
In addition, to ensure the integrability of $p(\mathbf{x},\mathbf{\alpha}),$ we enforce $\fa_5,\fa_9,\fa_{12},\fa_{14}$ and $\fa_{15},\ldots,\fa_{24}$ to be nonnegative. 

\subsection{Independent variables with known density}
For the first example, we chose $\fa_5=\fa_9=\fa_{12}=\fa_{14}=0.5,$ $\fa_{15}=\fa_{19}=\fa_{22}=\fa_{24}=1$ and the rest of coefficients were set equal to zero. With this choice $p(\mathbf{x})=\prod_{i=1}^4 \exp(-0.5x_i^2-x_i^4)/Z$ where $Z=\prod_{i=1}^4 \int_{-\infty}^{\infty}\exp(-0.5x_i^2-x_i^4)dx_i$ is the normalization constant. We used $N=10^6$ samples of this density to compute the moments  $T_i=\frac{1}{N} \sum_{j=1}^{N} \psi_k(\mathbf{x}_j) \; \text{for} \; k=1,\ldots,24.$ Before computing the moments, we normalized the samples by their empirical mean $\mu_i, \; i=1,\ldots,4,$ and standard deviation, $\sigma_i, \; i=1,\ldots,4,$ i.e. $x_{ij} \rightarrow (x_{ij}-\mu_i)/\sigma_i,$ for $j=1,\ldots,N.$ The goal was to estimate the parameters of an exponential density so that they reproduced the moments $T_i.$ For the naive global optimization, the generic particle filter optimization algorithm and the modified particle filter optimization algorithm we used $N_E=10^4$ samples to estimate all the necessary expectation values. The number of particles was set to $M=100.$ The initial condition for the $j$-th particle ($j=1,\ldots,M$) was chosen as 
\begin{gather*}
\fa^{(0)}_{j5}=.5(1-\eta_{j1}), \; \fa^{(0)}_{j9}=.5(1-\eta_{j2}), \\
\fa^{(0)}_{j12}=.5(1-\eta_{j3}), \; \fa^{(0)}_{j14}=.5(1-\eta_{j4}),\\
\fa^{(0)}_{j15}=1-.5\eta_{j5}, \; \fa^{(0)}_{j19}=1-.5\eta_{j6}, \\
\fa^{(0)}_{j22}=1-.5\eta_{j7}, \; \fa^{(0)}_{j24}=1-.5\eta_{j8},
\end{gather*} 
where $\eta_{j1},\ldots, \eta_{j8}$ are independent random variables uniformly distributed in $[0,1).$ The rest of the components of $\fa$ are set initially to zero. 

For the modified particle filter algorithm, we set the parameter $\gamma,$ which determines the value of $\lambda$ for the different particles in a batch, to the value $\gamma=1.1$ (see discussion at the end of Section \ref{application2.1}). We note that the same value of the parameter $\gamma$ was used for all four examples with equal success. The choice for $\gamma$ was not optimized, i.e., it was not found by trial and error. As stated before, it can be chosen based on knowledge of how the parameter $\lambda$ influences the behavior of the LM algorithm, in particular the determinant and condition number of the matrix $\tilde{J_j}^T \tilde{J_j}+ \lambda^{(t)}_j  diag (\tilde{J_j}^T \tilde{J_j}).$  

\begin{figure}
\centering
\epsfig{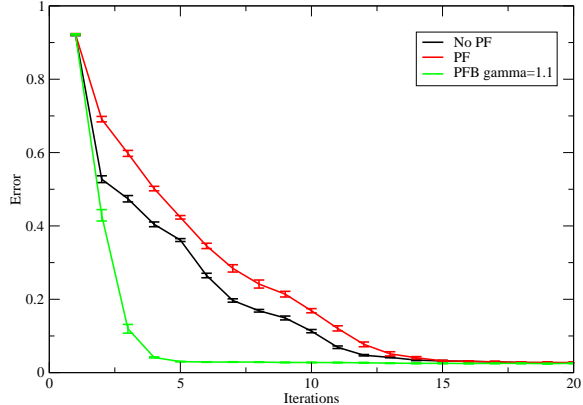}
\caption{First example: Independent variables with known density (see text for details). Evolution of the average error per moment as a function of LM iterations. No PF corresponds to the naive global optimization algorithm, PF corresponds to the generic particle filter optimization algorithm and PFB to the modified particle filter optimization algorithm.}
\label{plot_syn_diagonal}
\end{figure}

The most severe test for the particle filter and modified particle filter algorithms is to compare their error to the error of the naive global optimization algorithm when all three algorithms are started from identical initial conditions for the particles. In addition, to eliminate possible discrepancies in the algorithms' behavior due only to the variability inherent in Monte Carlo sampling we performed for each algorithm 10 different experiments, each one with 100 particles, with identical initial conditions and averaged the results over the 10 experiments. Figure \ref{plot_syn_diagonal} presents the evolution of the average error per moment $(\frac{2*\tilde{\epsilon}(\fa)}{24})^{1/2}$ as a function of LM iterations for the three algorithms. The error bars denote the standard deviation of the average over the 10 experiments. It is obvious from the error bars that there is not much variability between the different experiments. After a few iterations the standard deviation becomes an order of magnitude smaller than the average. 

While all three algorithms converge to approximately the same value for the error, which is practically equal to the Monte Carlo error, the speed at which they do so is very different. The naive algorithm converges faster than the generic particle filter algorithm. As we discussed before (end of Section \ref{application2}), this can happen. The modified particle filter algorithm significantly outperforms both the naive algorithm and the generic particle filter.  In particular, the modified algorithm has practically converged by iteration 5 while the other two algorithms need about 15 iterations to reduce the error to the same value. This is an increase in the convergence speed by more than $60\%.$ 

It is important to discuss the computational cost at which this convergence speedup is achieved. The difference between both versions of the particle filter algorithm and the naive algorithm is the addition of the filtering and resampling steps. However, this cost is negligible compared to the Monte Carlo sampling computational cost which is needed to setup the LM algorithm and compute the error at each step. For all the examples studied here, the computational cost to perform a fixed number of iterations with the modified particle filter algorithm is about $2\%$ more than the cost of the naive algorithm. Given the increase of the convergence speed this extra cost is well worthwhile. In particular, the extra cost of the modified particle filter algorithm is less than the cost of adding $2\%$ more particles in the naive algorithm. For this example this would mean adding 2 more particles to the naive algorithm.

\subsection{Dependent variables with known density}
\begin{figure}
\centering
\epsfig{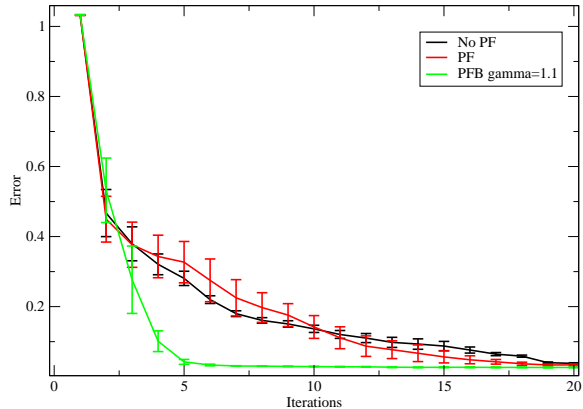}
\caption{Second example: Dependent variables with known density (see text for details). Evolution of the average error per moment as a function of LM iterations. No PF corresponds to the naive global optimization algorithm, PF corresponds to the generic particle filter optimization algorithm and PFB to the modified particle filter optimization algorithm.}
\label{plot_syn_full}
\end{figure}

For the second example, we chose $\fa_5=\fa_9=\fa_{12}=\fa_{14}=0.5,$ $\fa_{15}=\ldots=\fa_{24}=1$ and the rest of coefficients were set equal to zero. Thus, the random variables $x_1,\ldots,x_4$ are dependent. Figure \ref{plot_syn_full} presents the evolution of the average error per moment $(\frac{2*\tilde{\epsilon}(\fa)}{24})^{1/2}$ as a function of LM iterations for the three algorithms. The naive algorithm and the generic particle filter have comparable behavior while the modified particle algorithm outperforms both of them significantly. As in the first example, the error for the modified particle algorithm has practically converged by iteration 5 to a value comparable to the Monte Carlo sampling error while the other two algorithms need almost 20 iterations to do so. Thus, the convergence speedup of the modified particle algorithm is more than $70\%.$

\subsection{Sinusoidal signal with additive noise}
For the third and fourth examples we used the algorithms to estimate an exponential density that reproduces the moments $T_i=\frac{1}{N} \sum_{j=1}^{N} \psi_k(\mathbf{x}_j) \; \text{for} \; k=1,\ldots,24$ of an unknown density. The examples are motivated by training neural networks to represent time series \cite{bishop}. We will examine two cases of a signal corrupted by noise. In the third example we suppose that we have samples from a signal $u(t)=\sin(2\pi t)+\eta_t$ where $\eta_t$ is Gaussian white noise. At any instant $t,$ the noise $\eta_t \sim N(0,1).$ In the fourth example we suppose that we have samples from a signal $u(t)=\sin(2\pi (t+\eta))$ where $\eta \sim N(0,0.1),$ i.e. a signal with a random phase. For both cases we assume that the signal is given for $t \in [-1,1].$ We expand the signal in Legendre polynomials (which are orthogonal in [-1,1]) and keep only the first 4 terms in the expansion. Since the signal is random, the coefficients of the expansion are random. The random variables $x_1,\ldots,x_4$ are the coefficients of the first 4 Legendre polynomials. Exactly because the signal is a random function, we do not expect the expansion in Legendre polynomials to be an accurate one. The coefficients of the expansion are expected to fluctuate considerably from sample to sample. However, our purpose is to see how well an exponential density can represent the unknown density of the Legendre expansion coefficients. 

\begin{figure}
\centering
\epsfig{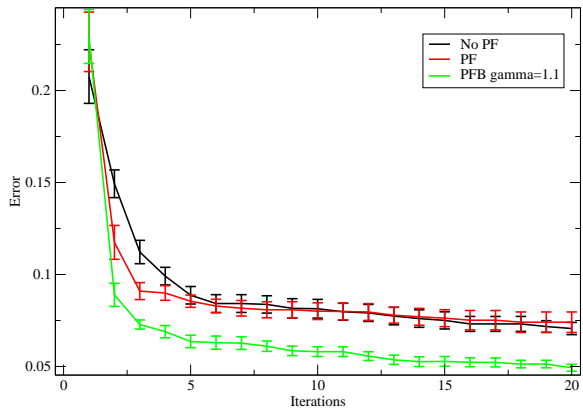}
\caption{Third example: Sinusoidal signal with additive noise (see text for details). Evolution of the average error per moment as a function of LM iterations. No PF corresponds to the naive global optimization algorithm, PF corresponds to the generic particle filter optimization algorithm and PFB to the modified particle filter optimization algorithm. }
\label{plot_gaussian}
\end{figure}

The parameters in the implementation of all three algorithms are the same as in the first two examples, except for the initial conditions. The initial condition for the $j$-th particle ($j=1,\ldots,M$) was chosen as 
\begin{gather*}
\fa^{(0)}_{j5}=1-\eta_{j1}, \; \fa^{(0)}_{j9}=1-\eta_{j2}, \\
\fa^{(0)}_{j12}=1-\eta_{j3}, \; \fa^{(0)}_{j14}=1-\eta_{j4},\\
\fa^{(0)}_{j15}=1-\eta_{j5}, \; \fa^{(0)}_{j19}=1-\eta_{j6}, \\
\fa^{(0)}_{j22}=1-\eta_{j7}, \; \fa^{(0)}_{j24}=1-\eta_{j8},
\end{gather*} 
where $\eta_{j1},\ldots, \eta_{j8}$ are independent random variables uniformly distributed in $[0,1).$ The rest of the components of $\fa$ are set initially to zero. 

Figure \ref{plot_gaussian} presents the evolution of the average error per moment $(\frac{2*\tilde{\epsilon}(\fa)}{24})^{1/2}$ as a function of LM iterations for the three algorithms. We see that the naive algorithm and the generic particle filter algorithm have comparable behavior. However, after 20 iterations, both of them have reduced the error to a value that is still about $30\%$ larger than the error of the modified particle filter algorithm. The modified particle filter algorithm reduces the error to about 1.5 times the Monte Carlo error.

\subsection{Sinusoidal signal with random phase}
As mentioned in the preceding section, we suppose that we have samples from a random phase signal $u(t)=\sin(2\pi (t+\eta))$ where $\eta \sim N(0,0.1),$ with $t \in [-1,1].$ As before, we expand the signal in Legendre polynomials and keep only the first 4 terms in the expansion. Since the signal is random, the coefficients of the expansion are random. The random variables $x_1,\ldots,x_4$ are the coefficients of the first 4 Legendre polynomials. The initial conditions for the particles are assigned in the same manner as in the third example.

The fourth example is quite more challenging than the third example. In particular, dependencies among the coefficients of the expansion, i.e. the random variables $x_1,\ldots,x_4,$ render the Jacobian $J$ singular. This means that the optimization problem admits nonzero solutions for the error vector $f_k(\fa)=E_{\mathbf{\alpha}}[\psi_k(\mathbf{x})]-\frac{1}{N} \sum_{j=1}^{N} \psi_k(\mathbf{x}_j), \quad k=1,\ldots,24.$ Even though the matrix $\tilde{J_j}^T \tilde{J_j}+ \lambda^{(t)}_j  diag (\tilde{J_j}^T \tilde{J_j})$ used in the calculation of the parameter vector increment is regularized, the best one can hope for is to keep reducing the error until the Jacobian becomes zero to within the arithmetic precision used ($10^{-16}$ in our case).

\begin{figure}
\centering
\epsfig{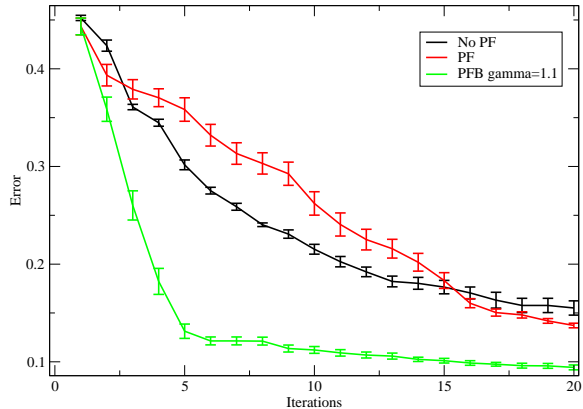}
\caption{Fourth example: Sinusoidal signal with random phase (see text for details). Evolution of the average error per moment as a function of LM iterations. No PF corresponds to the naive global optimization algorithm, PF corresponds to the generic particle filter optimization algorithm and PFB to the modified particle filter optimization algorithm.  }
\label{plot_nongaussian_ran_add}
\end{figure}

Figure \ref{plot_nongaussian_ran_add} presents the evolution of the average error per moment $(\frac{2*\tilde{\epsilon}(\fa)}{24})^{1/2}$ as a function of LM iterations for the three algorithms. The modified particle filter algorithm has practically converged by the 12th iteration with an error value that is about $30\%$ less than the error value achieved by the naive algorithm and the generic particle filter. However, the error value of the modified particle filter is still about 8 times larger than the Monte Carlo error. This is because, as we mentioned in the previous paragraph, this optimization problem admits nonzero solutions for the error vector $\mathbf{f}(\fa).$ Still, the convergence speedup of the modified particle filter algorithm for this example is about $40\%.$

\section{Discussion}\label{discussion}
We have presented a reformulation of stochastic global optimization as a filtering problem. In particular we have reformulated stochastic global optimization using a particle filter. This choice was based on the simplicity of implementation and flexibility of particle filters. We have exploited this flexibility to construct a modified particle filter filter that converges faster than naive global optimization. We have demonstrated the efficiency of the approach with several examples of varying difficulty. 

The flexibility allowed by the particle filter can be used to construct additional modified particle filters. For example, we can use ranking selection \cite{eiben} in the resampling step instead of proportionate selection. The advantage of ranking selection compared to proportionate selection is that it establishes a constant pressure of selecting the particle with the largest weight. This can be helpful when there exist several particles with almost equally large weights. In this case, proportionate selection may not be able (for a finite number of particles) to sample the particle with the largest weight. Another possible modification of the particle filter to enhance the exploration of the parameter space can come from the use of recombination procedures between particles after the resampling step has been performed. The motivation for such a procedure comes from genetic algorithms \cite{eiben}, where recombination constitutes probably the most important feature of such algorithms. However, there is no unique way to perform recombination and this can also be the weak point of genetic algorithms. Another modification of the particle filter that we have already mentioned is to allow the particles to evolve for a few iterations without enforcing the observation and resampling steps. This allows the particles to exhibit their potential as far as locating a local minimum is concerned. After the first few iterations we can start enforcing the observation and resampling steps to allocate more particles in the more promising areas of the parameter space. Such a modification can help avoid the potential problem of premature convergence (see also the discussion in Section \ref{application2}). Finally, we note that the particle filter algorithm can be used also to randomize a deterministic global optimization algorithm \cite{zhou}.

We hope that the algorithm proposed in this work will help in tackling the multitude of optimization problems originating from real-world applications.

\section*{Acknowledgements} 
I am grateful to Dr. V. Maroulas and Dr. J. Weare for many helpful discussions and comments and to M. Sendak for inspiration.

\end{document}